\documentclass[11pt]{amsart}
\usepackage{amssymb}
\usepackage{latexsym}
\usepackage{graphicx}
\usepackage{subfigure}
\usepackage{tikz}
\usepackage[linktocpage=true]{hyperref}
\usepackage{ulem}
\usepackage[left=3cm, right=3cm, top=3cm, bottom=3cm]{geometry}
\usepackage{color}

\hypersetup{ colorlinks   = true, 
urlcolor  = blue, 
linkcolor    = blue, 
citecolor   = blue 
}

\def\beq{\begin{equation}}
\def\eeq{\end{equation}}

\def\bm{\begin{matrix}}
\def\em{\end{matrix}}

\newcommand{\Z}{{\mathbb Z}}
\newcommand{\R}{{\mathbb R}}

\newcommand{\T}{{\mathbb T}}

\newcommand{\N}{{\mathbb N}}

\newcommand{\CC}{{\mathcal C}}

\newcommand{\CH}{{\mathcal H}}
\newcommand{\CI}{{\mathcal I}}

\newcommand{\CS}{{\mathcal S}}

\newcommand{\CZ}{{\mathcal Z}}

\newcommand{\CK}{{\mathcal K}}
\newcommand{\CO}{{\mathcal O}}

\newtheorem{Theorem}{Theorem}[section]
\newtheorem{remark}{Remark}[section]

\newtheorem{Lemma}{Lemma}[section]

\newtheorem{Proposition}{Proposition}[section]

\newcommand{\la}{\langle}
\newcommand{\ra}{\rangle}

\begin{document}

\title[]{Growth of Sobolev Norms in 1-d Quantum Harmonic Oscillator with Polynomial Time Quasi-periodic Perturbation}

\author{Jiawen Luo}
\address{School of Mathematical Sciences and Key Lab of Mathematics for Nonlinear Science, Fudan University, Shanghai 200433, China}
\email{20110180010@fudan.edu.cn}

\author{Zhenguo Liang}
\address{School of Mathematical Sciences and Key Lab of Mathematics for Nonlinear Science, Fudan University, Shanghai 200433, China}
\email{zgliang@fudan.edu.cn}
\thanks{Z. Liang was partially supported by National Natural Science Foundation of China (Grants No. 12071083)
and Natural Science Foundation of Shanghai (Grants No. 19ZR1402400).}

\author{Zhiyan Zhao}
\address{Universit\'e C\^ote d'Azur, CNRS, Laboratoire J. A. Dieudonn\'{e}, 06108 Nice, France}
\email{zhiyan.zhao@univ-cotedazur.fr}
\thanks{}

\begin{abstract}
We consider the one-dimensional quantum harmonic oscillator perturbed by a linear operator which is a polynomial of degree $2$ in $(x,-{\rm i}\partial_x)$, with coefficients quasi-periodically depending on time. By establishing the reducibility results, we
describe the growth of Sobolev norms. In particular, the $t^{2s}-$polynomial growth of ${\mathcal H}^s-$norm is observed in this model if the original time quasi-periodic equation is reduced to a constant Stark Hamiltonian.

\end{abstract}

%



\maketitle

\section{Introduction and main results}

The aim of this paper is to describe the growth of Sobolev norms for the time-dependent Schr\"odinger equation
\begin{equation}\label{eq_Schrodinger}
 {\rm i}\partial_t u=\frac{\nu(E)}{2}H_0u + W(E,\omega t, x, -{\rm i} \partial_x)u,\qquad x\in\R ,
\end{equation}
where, we assume that
\begin{itemize}
  \item for given $d\geq 1$, the frequencies $\omega\in \R^d$ satisfy the {\it Diophantine} condition (denoted by $\omega\in {\rm DC}_d(\gamma,\tau)$ for $\gamma>0$, $\tau>d-1$):
$$\inf_{j\in\Z}|\la n,\omega\ra-j|>\frac{\gamma}{|n|^\tau},\qquad \forall \ n\in\Z^d\setminus\{0\},$$
  \item the parameter $E\in {\CI}$, an interval $\subset \R$, and $\nu\in C^2({\CI},\R)$ satisfies
  $$|\nu'(E)|\geq l_1,\quad |\nu''(E)|\leq l_2,\qquad \forall \  E\in{\CI},$$
for some $l_1,l_2>0$,
  \item the operator $H_0$ is the {\it one-dimensional quantum harmonic oscillator}, i.e.
$$(H_0u)(x):=-(\partial_{x}^2 u)(x)+x^2\cdot u(x),\qquad \forall \ u\in L^2(\R),$$
  \item the perturbation $W(E,\theta, x,\xi)$ is a polynomial in $(x,\xi)$ of degree $2$, i.e.,
\begin{eqnarray}
W(E,\theta,x,\xi)&=& \frac12\big(a_{20}(E,\theta)x^2+2a_{11}(E,\theta)x\cdot\xi+a_{02}(E,\theta)\xi^2\big)\label{pertur_quadratic}\\
& & + \, b_1(E,\theta)x+b_2(E,\theta)\xi+c(E,\theta),\label{pertur_linear}
\end{eqnarray}
with the $\theta-$dependent coefficients $a_{20},\, a_{11}, \, a_{02}, \, b_1, \, b_2, \, c: {\CI}\times \T^d\to \R$ satisfying that
\begin{itemize}
  \item all above coefficients are real-analytic w.r.t. $\theta\in\T^d:=(\R/2\pi\Z)^d$, and can be extended to the complex neighborhood ${|\Im z|<r}$ (denoted by  $C_r^\omega(\T^d)$),
  \item the coefficients $a_{20}$, $a_{11}$, $a_{02}$ in quadratic part are $C^2$ w.r.t. $E\in{\CI}$ with $|\partial^m_E a_{20}(E,\cdot)|_r:=\sup_{|\Im z|<r}|\partial^m_E a_{20}(E,z)|$, $|\partial^m_E a_{11}(E,\cdot)|_r$, $|\partial^m_E a_{02}(E,\cdot)|_r$ small enough for every $E\in{\CI}$, for $m=0,1,2$.
  \end{itemize}
\end{itemize}

According to \cite{BGMR2018}, it is known that, for the polynomial of degree $2$, the classical-quantum correspondence is exact and it is possible to export techniques from classical to quantum mechanics.
With this idea, in \cite{LZZ2021}, for Eq. (\ref{eq_Schrodinger}) with a time quasi-periodic quadratic perturbation (i.e., $W = (\ref{pertur_quadratic})$), the authors showed that for a.e. $E\in {\CI}$ it is reducible, i.e., via a unitary transformation depending on time in an analytic quasi-periodic way, Eq.(\ref{eq_Schrodinger}) is conjugated to an equation independent of time. The idea of proof is based on the KAM theory of reducibility for quasi-periodic linear systems developed by Eliasson \cite{Eli1992}, in which the reducibility transformation is not necessarily close to identity and the type of normal form can be different from the original one.
In this paper, as a continuation of \cite{LZZ2021}, we consider the time quasi-periodic polynomial perturbation of degree $2$, i.e., $W = (\ref{pertur_quadratic})+ (\ref{pertur_linear})$.
Since the main difference from \cite{LZZ2021} is the presence of linear and constant terms in (\ref{pertur_linear}) in the perturbation $W$, it is sufficient to consider the reducibility for the time quasi-periodic lower-order part based on a constant quadratic part.

As an application of reducibility, the behaviors of Sobolev norms for the solutions to Eq. (\ref{eq_Schrodinger}) are described. In particular, the $t^{2s}$ polynomial growth of ${\CH}^s-$norm is possible for $E$'s in a zero-measure subset of ${\CI}$ if Eq.(\ref{eq_Schrodinger}) is reduced to a constant Stark Hamiltonian (p.530 of \cite{HLS86}). This is the remarkable difference from Theorem 2 of \cite{LZZ2021}.

More precisely, for $s\geq 0$, we consider the solution to Eq.(\ref{eq_Schrodinger}) in Sobolev space
$${\CH}^s:=\left\{\psi\in L^2(\R):H_0^{\frac{s}2}\psi \in L^2(\R)\right\}$$
and its Sobolev norm
$\|\psi\|_{s}:=\|H_0^{\frac{s}2} \psi\|_{L^2}$.
It is well known that, for $s\in \N$, the above definition of norm is equivalent to
\begin{equation}\label{norm_equiv}
\sum\limits_{\alpha+\beta\leq s\atop{\alpha,\beta\in\N} }\|x^{\alpha}\cdot\partial^{\beta} \psi\|_{L^2}.
\end{equation}
In view of Proposition 2.3 of \cite{BM2018} and Lemma 2.4 of \cite{YZ2004}, we get that, for a given $\psi\in {\CH}^s$,
\begin{equation}\label{norm_equiv-1}
\|\psi\|_{s}\simeq \|\psi\|_{H^s}+ \|x^{s} \psi\|_{L^2},
\end{equation}
where `` $H^s$ " means the standard Sobolev space and $\|\cdot\|_{H^s}$ is the corresponding norm. Hence, to calculate the norm $\|\psi\|_s$, $s\geq 0$, it is sufficient to focus on $\|x^{s} \psi\|_{L^2}$ for $s\geq0$ and $\|\psi^{(s)}\|_{L^2}$ for $s\in\N$.
In view of Theorem 1.2 of \cite{MR2017}, the solution to Eq. (\ref{eq_Schrodinger}) is globally well-posed in ${\CH}^s$.

\

Before the statement of our main results, let us review the previous works on the reducibility of harmonic oscillators and behaviors of solutions in Sobolev spaces.

The boundedness of Sobolev norms via reducibility, as well as the pure point nature of Floquet operator, was firstly considered for 1-d quantum harmonic oscillator (QHO for short) with time periodic smooth perturbations \cite{Com87, DLSV2002, EV83, Kuk1993}. As for time quasi-periodic perturbations, we can refer to \cite{GreTho11, Wang08, WLiang17} for bounded perturbations, and refer to \cite{Bam2018, Bam2017, BM2018, Liangluo19} for unbounded perturbations. In investigating the reducibility problems, KAM theory for 1-d PDEs has been well developed by Bambusi-Graffi \cite{BG2001} and Liu-Yuan \cite{LY2010}, which are contributed to KAM theory with unbounded perturbations.

Reducibility for PDEs in higher-dimensional case was initiated by Eliasson-Kuksin \cite{EK2009}. We can refer to \cite{GrePat16} and \cite{LiangWang19} for higher-dimensional QHO with bounded potential.
We also mention that some higher-dimensional results with unbounded perturbations have been recently obtained \cite{BGMR2018, BLM18, FGMP19, FG19, FGN19, Mon19}.
However, a general KAM theorem for higher-dimensional PDEs with unbounded perturbations is far from success.

As a supplement to the understanding of the long-time behavior of solutions in Sobolev spaces, it is worthy constructing unbounded solutions in Sobolev spaces, especially describing the growth rate with time of Sobolev norms.
For 1-d QHO with time periodic perturbation, the solutions with $t^{\frac{s}{2}}-$polynomial growth for ${\CH}^s-$norms are constructed by Delort \cite{Del2014}, Maspero \cite{Mas2018} and those with $t^{s}-$polynomial growth are constructed by Bambusi-Gr\'ebert-Maspero-Robert \cite{BGMR2018}.
In \cite{BGMR2018}, following Graffi-Yajima \cite{GY00}, the authors also constructed the solutions with $t^{s}-$polynomial growth for ${\CH}^s-$norms for the higher-dimensional QHO with time quasi-periodic perturbation which is linear in $x$ and $-{\rm i}\partial_x$ (via Theorem 3.3 of \cite{BGMR2018}).
For 1-d QHO with time quasi-periodic perturbation which is a quadratic form of $(x,-{\rm i}\partial_x)$, the $t^{s}-$polynomial growth and exponential growth for ${\CH}^s-$norms are shown by Liang-Zhao-Zhou \cite{LZZ2021}.
Recently, for 2-d QHO with perturbation which is decaying in $t$, Faou-Rapha\"el \cite{FR20} constructed a solution whose ${\CH}^1-$norm presents logarithmic growth with $t$. For 2-d QHO with perturbation being the projection onto Bargmann-Fock space, Thomann \cite{Thomann20} constructed explicitly a travelling wave whose Sobolev norm presents polynomial growth with $t$, based on the study in \cite{ST20} for linear Lowest Landau Level equations (LLL) with a time-dependent potential.
There are also many literatures, e.g., \cite{HM2020, Mas2021} which give growth of Sobolev norms for harmonic oscillators and anharmonic oscillators, and \cite{BGMR2019, BLM2021, BM2019, Bou99a, Bou99b, FZ12, MR2017, Wang2008}, which are relative to the upper growth bound of the solution in Sobolev space.

Compared with previous works, the novelty in this paper is the $t^{2s}-$polynomial growth for ${\CH}^s-$norms of solutions to Eq. (\ref{eq_Schrodinger}) because of the interaction between quadratic part and linear part in certain situations. The exponential growth is also shown in this paper as \cite{LZZ2021}, which gives another example for the optimal growth mentioned in Theorem 1.2 of \cite{MR2017}.

\subsection{Growth of Sobolev norms}

The main result of this paper is about the growth of Sobolev norms of solutions to Eq. (\ref{eq_Schrodinger}).

\begin{Theorem}\label{thm_Schro_sobolev}
There exists $\epsilon_*=\epsilon_*(\gamma,\tau,r,d,l_1,l_2)>0$ such that if
$$\sup_{E\in\mathcal{I}}\max_{m=0,1,2}\{|\partial_E^m a_{20}|_r, \, |\partial_E^m a_{11}|_r, \, |\partial_E^m a_{02}|_r\}=:\epsilon_0\le\epsilon_*,$$
then there exists a subset
$$\mathcal{O}_{\epsilon_0}=\bigcup_{j\in\mathbb{N}}\Lambda_j\subset\overline{\mathcal{I}}$$
with $\Lambda_j$ being closed intervals \footnote{The ``closed interval" here is interpreted in a more general sense, i.e., it can be degenerated to a point instead of a positive-measure subset of $\R$.} and $Leb(\mathcal{O}_{\epsilon_0}) < \epsilon_0^{\frac{1}{40}}$, such that the following holds for the solution $u(t)$ to Eq. (\ref{eq_Schrodinger}) for $t\geq0$ and for non-vanishing $u(0)\in {\CH}^s$, $s\geq 0$.

There exist two constants $c, \, C>0$, depending on $s$, $E$ and $u(0)$, such that
\begin{enumerate}
\item For a.e. $E\in{\CI}\setminus{\CO}_{\epsilon_0}$, $c \leq \|u(t)\|_{s}\leq C$.
\item If ${\rm Leb}(\Lambda_j)>0$, then
\begin{itemize}
  \item for $E\in{\rm int}\Lambda_j$, there exists $\lambda=\lambda_E >0$ such that
  \begin{equation}\label{exponential_growth}
  c  e^{\lambda st}  \leq \|u(t)\|_{s} \leq C  e^{\lambda st},
  \end{equation}
  \item for $E\in\partial \Lambda_j\setminus \partial{\CI}$, there exist $\kappa=\kappa_E \in\R\setminus\{0\}$ and $\iota=\iota_E\in\R$ such that
   \begin{equation}\label{polynomial_growth}
  c( |\kappa|^st^s+|\iota\kappa|^{s}t^{2s})\le\|u(t)\|_s\le C(1+|\kappa|^st^s+|\iota\kappa|^{s}t^{2s}).
  \end{equation}
\end{itemize}
If ${\rm Leb}(\Lambda_j)=0$, then for $E\in \Lambda_j$, there exists $\iota=\iota_E\in\R$ such that
\begin{equation}\label{polynomial_growth-1}
c |\iota|^{s}t^{s}   \leq \|u(t)\|_{s} \leq C(1+|\iota|^{s}t^{s}).
  \end{equation}
\end{enumerate}
\end{Theorem}

\begin{remark}\label{rem_iota}
According to (\ref{polynomial_growth}), we see that the $t^{2s}-$polynomial growth of Sobolev norm occurs when the constants $\kappa$ and $\iota$ in (\ref{polynomial_growth}) do not vanish, which means that Eq. (\ref{eq_Schrodinger}) is reduced to a constant Stark Hamiltonian (see Theorem \ref{thm_Schro_reduc} and (\ref{stark}) in Section \ref{sec_proof_reduc}).
In Section 1.3 of \cite{LZZ2021}, there are several concrete examples of time quasi-periodic quadratic perturbations for which $\kappa$, as well as $\lambda$ in (\ref{exponential_growth}), does not vanish. On the other hand, we can see from the proof that the non-vanishing time quasi-periodic linear terms $b_1(E,\theta)x$ and $b_2(E,\theta)\xi$ in $W(E,\theta,x,\xi)$, for which $\iota$ does not vanish, universally exist.
\end{remark}

Let us compare the behaviors of solutions in Theorem \ref{thm_Schro_sobolev} with those in Theorem 2 of \cite{LZZ2021}, since the present paper is a continuation of \cite{LZZ2021}.
\begin{itemize}
  \item In Theorem \ref{thm_Schro_sobolev}, the perturbation $W=(\ref{pertur_quadratic})+(\ref{pertur_linear})$ in Eq. (\ref{eq_Schrodinger}) is a polynomial of $(x,-{\rm i}\partial_x)$ of degree $2$.
  \item In Theorem 2 of \cite{LZZ2021}, the perturbation $W=(\ref{pertur_quadratic})$ in the equation is a quadratic form of $(x,-{\rm i}\partial_x)$.
\end{itemize}
The classification of parameters $E$'s in Theorem \ref{thm_Schro_sobolev} is the same with that of \cite{LZZ2021}, which is based on the reducibility and rotation number of quasi-periodic linear systems. Some more descriptions are given in Section \ref{redu_qp_lin_sys}.

\noindent
\begin{tabular}{|l|l|l|l|}
\hline
&  $W=$ (\ref{pertur_quadratic}) as in \cite{LZZ2021} & $W=$ (\ref{pertur_quadratic}) + (\ref{pertur_linear}) as in Theorem \ref{thm_Schro_sobolev} \\[2mm]
 \hline
for a.e. $E\in {\CI}\setminus {\CO}_{\epsilon_0}$  & boundedness & boundedness \\[2mm]
\hline
for $E\in{\rm int}\Lambda_j$ if ${\rm Leb}(\Lambda_j)>0$  & exponential growth & exponential growth \\[2mm]
\hline
for $E\in\partial \Lambda_j\setminus \partial{\CI}$  if ${\rm Leb}(\Lambda_j)>0$   & $t^s-$growth & possible $t^{2s}-$growth \\[2mm]
\hline
for $E\in\Lambda_j$ if ${\rm Leb}(\Lambda_j)=0$     & boundedness & possible $t^{s}-$growth \\[2mm]
\hline
\end{tabular}
We see that, the main differences are the possible $t^{2s}-$polynomial growth for the ${\CH}^s-$norms in (\ref{polynomial_growth}) for $E\in\partial \Lambda_j\setminus \partial{\CI}$ when ${\rm Leb}(\Lambda_j)>0$ and the possible $t^{s}-$polynomial growth in (\ref{polynomial_growth-1}) for $E\in \Lambda_j$ when $\Lambda_j$ is degenerated to a point.

Theorem \ref{thm_Schro_sobolev} indeed  implies Theorem 2 of \cite{LZZ2021}. In Theorem \ref{thm_Schro_sobolev}, the possible $t^{2s}-$growth shown in (\ref{polynomial_growth}), as well as the possible $t^{s}-$growth shown in (\ref{polynomial_growth-1}), corresponds to a normal form with a linear (w.r.t. $x$ or $-{\rm i}\partial_x$) term with coefficient $\iota$. The growth in (\ref{polynomial_growth}) (resp. in (\ref{polynomial_growth-1})) is a real $t^{2s}-$growth (resp. $t^{s}-$growth) only if the coefficient $\iota\neq 0$, since $|\iota|$ is a factor of coefficient of $t^{2s}$ in (\ref{polynomial_growth}) (resp. $t^{s}$ in (\ref{polynomial_growth-1})).
In particular, if the perturbation is quadratic as in \cite{LZZ2021}, i.e., the linear terms in (\ref{pertur_linear}) vanish, then $\iota$ also vanishes and hence the growth in (\ref{polynomial_growth}) (resp. in (\ref{polynomial_growth-1})) in Theorem 1.1 is exactly the $t^s-$growth (resp. the boundedness) given in Theorem 2 of \cite{LZZ2021} for $E\in\partial\Lambda_j\setminus \partial{\CI}$ (resp. for $E\in\Lambda_j$ if $\Lambda_j$ is degenerated to a point).

It was well believed that the quantum Hamiltonian (\ref{eq_Schrodinger}), perturbed by a polynomial of $(x,-{\rm i}\partial_x)$ of degree $2$, was characterized by the quadratic terms in the perturbation and the linear terms would not affect the property of the systems in front of quadratic terms. Indeed, according to the above table, it is true for a.e. parameter $E\in {\CI}$.
However, for countably infinite $E$'s, the behavior of solutions in the Sobolev space could be changed because of the interaction between the quadratic terms and the linear terms.

%
%

\subsection{Time evolution for $1-$d Stark operator}

On route to the proof of $t^{2s}-$growth in (\ref{polynomial_growth}), we discover the dynamics of the time evolution for one-dimensional Stark operator concerning the growth of the weighted $L^2-$norm of its time evolution.
More precisely, for the $1-$d Stark operator ${\CS}_a$, $a\in\R\setminus\{0\}$, on $L^2(\R)$, defined as
$$({\CS}_au)(x)=-u''(x)+ a x u(x),\quad x\in\R ,$$
we consider the weighted $L^2-$norm of the time evolution
$$\left(\int_{\R} x^{2s} \left|\left(e^{-{\rm i}t{\CS}_a}u_0\right)(x)\right|^2 dx\right)^{\frac12},\quad u_0\in{\CH}^s,\quad s\geq 0.$$

\begin{Theorem}\label{super-ballistic}
For any non-vanishing $u_0\in{\CH}^s$, we have that
$$\lim_{t\to\infty}\frac{1}{t^{2s}}\left(\int_{\R} x^{2s} \left|\left(e^{-{\rm i}t{\CS}_a}u_0\right)(x)\right|^2 dx\right)^{\frac12}=  |a|^s \|u_0\|_{L^2}.$$
\end{Theorem}

\begin{remark} It has been well conjectured that absolute continuity implies ballistic motion (see, for example, \cite{AK1998}), i.e., the above weighted $L^2-$norm of the time evolution presents the $t^s-$growth with $t$. Since the $1-$d Stark operator has purely absolutely continuous spectrum on $\R$ (see e.g., \cite{AvHe1977, Kis1999}), Theorem \ref{super-ballistic} provides a counter-example for this conjecture.
\end{remark}

\subsection{Example with the $t^{2s}$ growth of $\mathcal{H}^s$ norm}

In the following, we will present a concrete time periodic quantum harmonic oscillator as an example, whose propagator has the $t^{2s}-$growth of $\mathcal{H}^s-$norm.

Consider the time periodic Schr\"odinger equation
\begin{equation}\label{exameq1}
{\rm i}\partial_t u =\frac{1}{2}H_0u+ W(t,x,-{\rm i}\partial_x)u, \quad u(0)\in\mathcal{H}^s,
\end{equation}
where $H_0$ is the $1-$d QHO as in  (\ref{eq_Schrodinger}) and $W(t,x,\xi)$ is a polynomial of degree $2$ about $(x,\xi)$ and periodic in $t$ s.t.
\begin{eqnarray*}
W(t,x,-{\rm i}\partial_x)&=&\frac\kappa2\left(-\cos^2(t)\cdot x^2-\cos(t) \sin(t) (x\cdot{\rm i}\partial_x+{\rm i}\partial_x\cdot x)+\sin^2(t)\cdot \partial^2_x\right)\\
&& + \, 2{\rm i}\iota\cos(t)\cdot \partial_x,
\end{eqnarray*}
for real constants $\kappa,\iota\ne0$.

\begin{Theorem}\label{examthm}
For $s\geq 0$ and non-vanishing $u(0)\in\mathcal{H}^s$, there are two constants $c$, $C >0$, depending on $s$ and $u(0)$, such that the solution to Eq. (\ref{exameq1}) fulfills
$$c |\kappa\iota|^s t^{2s}\le\|u(t)\|_{s}\le C(1+|\kappa\iota|^s t^{2s}),\quad\forall \ t\geq 0.$$
\end{Theorem}

\

 The rest of paper will be organized as follows. In Section \ref{sec_prelimi}, which serves as a preliminary section, we recall some known results on the relation between classical Hamiltonian and quantum Hamiltonian, and on the reducibility of the quasi-periodic linear systems. We consider the reducibility of Eq. (\ref{eq_Schrodinger}) in Section \ref{sec_reducibility} by the reducibility of the corresponding classical Hamiltonians and affine systems. In Section \ref{sec_Growth}, according to the type of reduced quantum Hamiltonians, we calculate the growth of Sobolev norms of solutions to Eq. (\ref{eq_Schrodinger}), and give the proof of Theorem \ref{thm_Schro_sobolev} and \ref{super-ballistic}.
In Section \ref{sec_pr_examples}, the proof of Theorem \ref{examthm} is given by a reducibility to the Stark Hamiltonian for the time periodic Schr\"odinger equation (\ref{exameq1}).

\

\noindent {\bf Acknowledgements.}
Z. Zhao would like to thank the Key Lab of Mathematics for Nonlinear Science of Fudan University for its hospitality during his visits in 2021.

 \section{Preliminaries}\label{sec_prelimi}

\subsection{Classical and quantum Hamiltonians}
Let us recall the definition of Weyl quantization, which relates the classical and quantum mechanics. The conclusions listed in this section can also be found in \cite{BGMR2018}.

The Weyl quantization is the operator ${\rm Op}^W:f\mapsto f^W$ for any symbol $f=f(x,\xi)$, with $x,\xi\in\R^n$, $n\geq 1$, where $f^{W}$ is the Weyl operator of $f$:
$$\left(f^{W} u\right)(x)=\frac{1}{(2\pi)^n}\int_{y, \, \xi\in\R^n} e^{{\rm i}(x-y)\xi} f\left(\frac{x+y}{2},\xi\right) u(y) \, dy  \, d\xi,\qquad \forall \ u\in L^2(\R^n).$$
In particular, if $f$ is a polynomial of degree at most $2$ in $(x,\xi)$, then $f^W$ is a polynomial of degree at most $2$ in $(x,-{\rm i}\partial_x)$ after the symmetrization.

For the $1-$parameter family of Hamiltonian $\chi(t, x, \xi )$, with $t$ an external parameter, let $\phi^\tau(t,x,\xi)$ be the time $\tau-$flow it generates, precisely the
solution to
$$\frac{dx}{d\tau}=\frac{\partial\chi}{\partial\xi}(t, x, \xi ),\qquad \frac{d\xi}{d\tau}=-\frac{\partial\chi}{\partial x}(t, x, \xi).$$
The time-dependent coordinate transformation
\begin{equation}\label{time1}
(x,\xi)=\phi^1\left(t,\tilde x,{\tilde\xi}\right)=\left.\phi^{\tau}\left(t,\tilde x,{\tilde\xi}\right)\right|_{\tau=1}
\end{equation}
transforms a Hamiltonian system with
Hamiltonian $h$ into another one with Hamiltonian $g$ given by
$$g\left(t,\tilde x,\tilde\xi\right)=h\left(\phi^1\left(t,\tilde x,\tilde\xi\right)\right)-\int_0^1 \frac{\partial\chi}{\partial t}\left(t,\phi^{\tau}\left(t,\tilde x,\tilde\xi\right) \right) d\tau.$$

\begin{Lemma}
Let $\chi(t, x, \xi)$ be a polynomial of degree at most $2$ in $(x,\xi)$ with smooth time-dependent coefficients.
If the Weyl operator $\chi^W(t, x, -{\rm i}\partial_x)$ is self-adjoint for any fixed $t$, then the transformation
\begin{equation}\label{tran}
\psi=e^{- {\rm i}\chi^W(t, x, -{\rm i}\partial_x)}\tilde\psi
\end{equation}
transforms the equation ${\rm i}\partial_t\psi=H\psi$ into ${\rm i}\partial_t\tilde\psi=G\tilde\psi$ with
\begin{eqnarray*}
G&:=&e^{{\rm i}\chi^W(t, x, -{\rm i}\partial_x)}He^{-{\rm i}\chi^W(t, x, -{\rm i}\partial_x)}\\
& & - \, \int_0^1 e^{{\rm i}\tau\chi^W(t, x, -{\rm i}\partial_x)}\left(\partial_t \chi^W(t, x, -{\rm i}\partial_x)\right)e^{-{\rm i}\tau\chi^W(t, x, -{\rm i}\partial_x)}d\tau.
\end{eqnarray*}
Moreover, if the transformation (\ref{time1}) transforms a classical system with Hamiltonian $h$ into
a system with Hamiltonian $g$, then the transformation (\ref{tran}) transforms the quantum Hamiltonian system $h^W$ into $g^W$.
\end{Lemma}


\begin{Proposition}\label{Prop_hami}
Let $\chi(\theta,x,\xi)$ be a polynomial of degree at most $2$ in $(x,\xi)$ with real coefficients depending in a $C^\infty-$way on $\theta\in \T^d$. For every $\theta\in \T^d$, the Weyl operator $\chi^W(\theta,x, -{\rm i}\partial_x)$ is self-adjoint in $L^2(\R^n)$ and $e^{-{\rm i}\tau\chi^W(\theta,x, -{\rm i}\partial_x)}$ is unitary in $L^2(\R^n)$ for every $\tau\in\R$.
Furthermore, for any $s\geq 0$, there exist $c'$, $C' > 0$ depending on $\|[H_0^s,\chi^W(\theta,x, -{\rm i}\partial_x)]H_0^{-s}\|_{L^2\mapsto L^2}$ and $s$, such that
   \begin{equation}\label{norm_U}
      c'\|\psi\|_{s}\leq \|e^{-{\rm i}\tau\chi^W(\theta,x,-{\rm i}\partial_x)}\psi\|_{s}\leq C'\|\psi\|_{s},\qquad \tau\in [0,1], \quad \theta\in\T^d.
    \end{equation}
\end{Proposition}

%
%
%


\subsection{Reducibility of quasi-periodic linear systems}\label{redu_qp_lin_sys}

We recall the reducibility of the ${\rm sl}(2,\R)$ quasi-periodic linear system $(\omega, A_E+F_E(\cdot))$:
\begin{equation}\label{linear_system_pr}
\frac{d}{dt}\left(\begin{array}{c}
x \\[1mm] \xi \end{array}\right)=(A_E+F_E(\omega t))\left(\begin{array}{c}
x \\[1mm] \xi \end{array}\right),
\end{equation}
where, with $\nu(E)$, $a_{20}$, $a_{11}$, $a_{02}$ as in Eq. (\ref{eq_Schrodinger}), the matrices are defined as
\begin{itemize}
\item $A_E :=\left(\begin{array}{cc}
0 & \nu(E) \\[1mm] -\nu(E) & 0 \end{array}\right)\in {\rm sl}(2,\R)$,
\item $ F_E(\cdot):=\left(\begin{array}{cc}
a_{11}(E,\cdot) & a_{02}(E,\cdot) \\[1mm] -a_{20}(E,\cdot) & -a_{11}(E,\cdot) \end{array}\right)\in C_r^{\omega}(\T^d, {\rm sl}(2,\R))$.
\end{itemize}
The reducibility of the above linear system was proved by Eliasson \cite{Eli1992}. We summarize the needed results in the following proposition as \cite{LZZ2021}.

Before stating the precise result on the reducibility of the above linear system,
we introduce the concept of rotation number. The {\it rotation number}
of quasi-periodic ${\rm sl}(2,\R)-$linear system $(\omega, A_E+F_E(\cdot))$
 is defined as
$$\rho_E:={\rm rot}(\omega, \, A_E+F_E(\cdot))=\lim_{t\to\infty}\frac{\arg(\Phi_E^t X)}{t},\quad \forall \ X\in \mathbb{R}^2\setminus\{0\},$$
where $\Phi_E^t$ is the basic
matrix solution and $\arg$ denotes the angle. The rotation number
$\rho_E$ is well-defined and it does not depend on $X$ \cite{JM82}.

\begin{Proposition}\label{prop_eliasson} [Proposition 5 of \cite{LZZ2021}] There exists $\epsilon_*=\epsilon_*(r,\gamma,\tau,d,l_1,l_2)>0$ such that if
\begin{equation}\label{small_F_0}
\sup_{E\in\mathcal{I}}\max_{m=0,1,2}\{|\partial_E^m a_{20}|_r, \, |\partial_E^m a_{11}|_r, \, |\partial_E^m a_{02}|_r\}=:\epsilon_0\le\epsilon_* ,
\end{equation}
then the following holds for the quasi-periodic linear system $(\omega, \, A_E+F_E(\cdot))$.
\begin{enumerate}
  \item [(1)] For a.e. $E\in{\CI}$, $(\omega, \, A_E+F_E(\cdot))$ is reducible. More precisely,
 there exist $B_E\in{\rm sl}(2,\R)$ and  $Z^E_j\in C^\omega(2\T^d, {\rm sl}(2,\R))$, $j=0,1,\cdots,K_E$, such that
 \begin{equation}\label{reducibility_sl2R}
   \frac{d}{dt}\left(\prod_{j=0}^{K_E} e^{Z^E_j(\omega t)}\right)=\left(A_E+F_E(\omega t)\right)\left(\prod_{j=0}^{K_E} e^{Z^E_j(\omega t)}\right)-\left(\prod_{j=0}^{K_E} e^{Z^E_j(\omega t)}\right)B_E.
 \end{equation}
  \item [(2)] The rotation number $\rho=\rho_E$ is monotonic on $\CI$. For any $k\in\Z^d$,
  $$\tilde\Lambda_k:=\left\{E\in\overline{\CI}:\rho_E=\frac{\la k,\omega\ra}{2}\right\}  \   \footnote{$\tilde\Lambda_k$ is defined as an empty set if the closed interval $\rho^{-1}\left(\frac{\la k,\omega\ra}{2}\right)$ does not intersect $\overline{\CI}$.} $$ is a closed interval, and we have
\begin{equation}\label{measure_esti}
\sum_{k\in\Z^d}{\rm Leb}(\tilde\Lambda_k)<\epsilon_0^{\frac{1}{40}}.
\end{equation}
  \item [(3)] For every $E\in \tilde\Lambda_k=:[a_k,b_k]$, $(\omega, \, A_E+F_E(\cdot))$ is reducible and the matrix $B_E\in {\rm sl}(2,\R)$ in (\ref{reducibility_sl2R}) satisfies
 \begin{itemize}
   \item if $a_k=b_k$, then $B_E=\left(\begin{array}{cc}
                                       0 & 0 \\
                                       0 & 0
                                     \end{array}\right)$;
   \item if $a_k<b_k$, then
\begin{itemize}
  \item ${\rm det} B_E<0$ for $E\in(a_k,b_k)$,
  \item ${\rm det} B_E=0$ for $E=a_k,b_k$ and $E\not\in\partial {\CI}$.
\end{itemize}
 \end{itemize}
   \item [(4)] For a.e. $E\in{\CI}\setminus\bigcup_k \tilde\Lambda_k$, $(\omega, \, A_E+F_E(\cdot))$ is reducible and the matrix $B_E\in {\rm sl}(2,\R)$ in (\ref{reducibility_sl2R}) satisfies ${\rm det} B_E>0$.
\end{enumerate}
\end{Proposition}

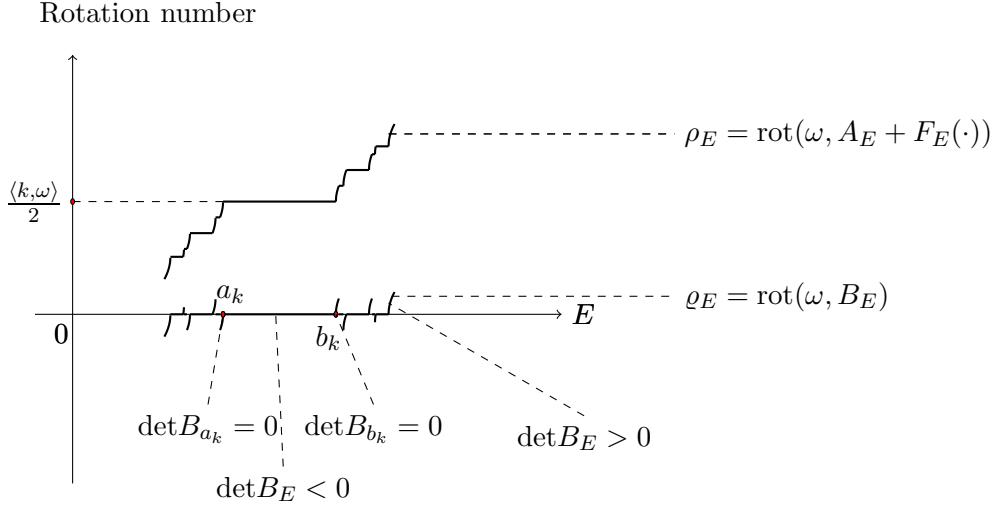
\begin{figure}\label{figure_rot_num}
\begin{center}
\begin{tikzpicture}[yscale=1.5]
\draw [->] (-0.5,0) -- (6.5,0);
\draw [->] (0,-1.5) -- (0,2.3);
\draw [thick,domain=1.22:1.3] plot (\x, {-0.7*sqrt(1.3-\x)});
\draw [thick,domain=1.3:1.47] plot (\x, {0});
\draw [thick,domain=1.47:1.48] plot (\x, {0.07-0.7*sqrt(1.48-\x)});
\draw [thick,domain=1.48:1.52] plot (\x, {0});
\draw [thick,domain=1.52:1.56] plot (\x, {-0.7*sqrt(1.56-\x)});
\draw [thick,domain=1.56:1.86] plot (\x, {0});
\draw [thick,domain=1.86:1.9] plot (\x, {0.14-0.7*sqrt(1.9-\x)});
\draw [thick,domain=1.9:1.96] plot (\x, {0});
\draw [thick,domain=1.96:2] plot (\x, {-0.7*sqrt(2-\x)});
\draw [thick,domain=2:3.5] plot (\x, {0});
\draw [thick,domain=3.5:3.54] plot (\x, {0.7*sqrt(\x-3.5)});
\draw [thick,domain=3.54:3.6] plot (\x, {0});
\draw [thick,domain=3.6:3.64] plot (\x, {-0.14+0.7*sqrt(\x-3.6)});
\draw [thick,domain=3.64:3.94] plot (\x, {0});
\draw [thick,domain=3.94:3.98] plot (\x, {0.7*sqrt(\x-3.94)});
\draw [thick,domain=3.98:4.02] plot (\x, {0});
\draw [thick,domain=4.02:4.03] plot (\x, {-0.07+0.7*sqrt(\x-4.02)});
\draw [thick,domain=4.03:4.2] plot (\x, {0});
\draw [thick,domain=4.2:4.28] plot (\x, {0.7*sqrt(\x-4.2)});
\draw [fill=red] (3.5,0) circle [radius=0.025];
\draw [fill=red] (2,0) circle [radius=0.025];
\node [below] at (-0.15,0) {$0$};
\node [right] at (6.5,0) {$E$};
\draw [thick,domain=1.22:1.3] plot (\x, {-0.49+1-0.7*sqrt(1.3-\x)});
\draw [thick,domain=1.3:1.47] plot (\x, {-0.49+1});
\draw [thick,domain=1.47:1.48] plot (\x, {-0.42+1-0.7*sqrt(1.48-\x)});
\draw [thick,domain=1.48:1.52] plot (\x, {-0.42+1});
\draw [thick,domain=1.52:1.56] plot (\x, {-0.28+1-0.7*sqrt(1.56-\x)});
\draw [thick,domain=1.56:1.86] plot (\x, {-0.28+1});
\draw [thick,domain=1.86:1.9] plot (\x, {-0.14+1-0.7*sqrt(1.9-\x)});
\draw [thick,domain=1.9:1.96] plot (\x, {-0.14+1});
\draw [thick,domain=1.96:2] plot (\x, {-0.7*sqrt(2-\x)+1});
\draw [thick,domain=2:3.5] plot (\x, {0+1});
\draw [dashed] (0,1)--(2,1);
\draw [thick,domain=3.5:3.54] plot (\x, {0.7*sqrt(\x-3.5)+1});
\draw [thick,domain=3.54:3.6] plot (\x, {0.14+1});
\draw [thick,domain=3.6:3.64] plot (\x, {0.14+0.7*sqrt(\x-3.6)+1});
\draw [thick,domain=3.64:3.94] plot (\x, {0.28+1});
\draw [thick,domain=3.94:3.98] plot (\x, {0.28+1+0.7*sqrt(\x-3.94)});
\draw [thick,domain=3.98:4.02] plot (\x, {0.42+1});
\draw [thick,domain=4.02:4.03] plot (\x, {0.42+1+0.7*sqrt(\x-4.02)});
\draw [thick,domain=4.03:4.2] plot (\x, {0.49+1});
\draw [thick,domain=4.2:4.28] plot (\x, {0.49+1+0.7*sqrt(\x-4.2)});
\draw [dashed] (4.2,1.6)--(8,1.6);
\draw [fill=red] (3.5,0) circle [radius=0.025];
\draw [fill=red] (2,0) circle [radius=0.025];
\draw [fill=red] (0,1) circle [radius=0.025];
\draw [dashed] (4.24,0.16)--(8,0.16);
\node [right] at (8,1.6) {$\rho_E={\rm rot}(\omega,A_E+F_E(\cdot))$};
\node [right] at (8,0.16) {$\varrho_E={\rm rot}(\omega,B_E)$};
\node [above] at (1,2.5) {Rotation number};
\node [left] at (0,1) {$\frac{\la k, \omega \ra}{2}$};
\node [below] at (-0.15,0) {$0$};
\node [right] at (6.5,0) {$E$};
\draw [dashed] (4.2,1.6)--(8,1.6);
\draw [dashed] (4,-0.8)--(3.5,0);
\node [below] at (3.4,-0.01) {$b_k$};
\draw [dashed] (1.8,-0.8)--(2,0);
\node [above] at (2.1,0.01) {$a_k$};
\draw [dashed] (2.8,-1.35)--(2.7,0);
\node [below] at (2.8,-1.35) {${\rm det}B_E<0$};
\node [below] at (1.8,-0.8) {${\rm det}B_{a_k}=0$};
\node [below] at (4,-0.8) {${\rm det}B_{b_k}=0$};
\draw [dashed] (6.8,-0.9)--(4.2,0.1);
\node [below] at (6.8,-0.9) {${\rm det}B_E>0$};
\end{tikzpicture}
\caption{Rotation numbers of linear systems $(\omega, A_E+F_E(\cdot))$ and $(\omega, B_E)$}
\label{f.graph}
\end{center}
\end{figure}

\begin{remark}\label{remark_rot_num} The original statement about reducibility in \cite{Eli1992} is that the quasi-periodic linear system $(\omega, A_E+F_E(\cdot))$ is reducible if the rotation number $\rho_E ={\rm rot}(\omega, \, A_E+F_E(\cdot))$ is Diophantine w.r.t. $\omega$, i.e., there exist $\sigma'$, ${\CK}'>0$, depending on $E$, such that
\begin{equation}\label{Dio_rho}
\left|\rho_E -\frac{\la k,\omega\ra}{2}\right|\geq \frac{{\CK}'}{|k|^{\sigma'}},\quad \forall \ k\in\Z^d\setminus\{0\},
\end{equation}
or $\rho_E $ is resonant, i.e., there exists $k\in\Z^d$ such that
 $$\rho_E =\frac{\la k,\omega\ra}{2}.$$
 According to \cite{Eli1992} and \cite{JM82}, the rotation number $\rho_E$ of the original system $(\omega, A_E+F_E(\cdot))$ is continuous and non-decreasing on ${\CI}$ because of the regularity and monotonicity of $\nu(E)$.
As shown symbolically in the graph locally around the interval $\tilde\Lambda_k=[a_k, b_k]$, $\rho_E$ is usually a Cantor staircase function of $E$.
However, the rotation number of the reduced system $\varrho_E :={\rm rot}(\omega, B_E)$ has no such continuity as $\rho_E$ since a renormalisation in the reducibility transformation makes $\varrho_E$ vanishing at every $\tilde\Lambda_k$. Indeed, for every $E\in{\CI}$, $\varrho_E ={\Im}\sqrt{{\rm det}B_E}$, and there is some $k_*=k_*(E)\in\Z^d$, constant at every $\tilde\Lambda_k$, such that
$$\varrho_E =\rho_E -\frac{\la k_*,\omega \ra}{2}.$$
Hence, if $\rho_E$ is non-vanishing and satisfies the Diophantine condition (\ref{Dio_rho}), then there exist $\sigma$, ${\CK}>0$, depending on $E$, such that
 \begin{equation}\label{Dio_varrho}
\left|\varrho_E-\frac{\la k,\omega\ra}{2}\right|=\left|\rho_E-\frac{\la k_*+ k,\omega\ra}{2}\right|\geq \frac{{\CK}}{1+|k|^{\sigma}},\quad \forall \ k\in\Z^d.
\end{equation}

\end{remark}

Let $\{\Lambda_j\}_{j\in\N}$ be the intervals $\tilde\Lambda_k$'s intersecting $\overline{\CI}$ and let
$${\CO}_{\epsilon_0}:=\bigcup_{j\in\N} \Lambda_j =\bigcup_{k\in\Z^d} \tilde\Lambda_k.$$
These are indeed the subsets stated in Theorem \ref{thm_Schro_sobolev}, as in Theorem 1 and 2 of \cite{LZZ2021}.

\section{Reducibility in classical and quantum Hamiltonians}\label{sec_reducibility}

In this section, we study the reducibility of the quantum Hamiltonian (\ref{eq_Schrodinger}) for a.e. $E\in{\CI}$ via the corresponding classical Hamiltonian. We have the following statement.

\begin{Theorem}\label{thm_Schro_reduc} Consider the time-dependent Schr\"odinger equation (\ref{eq_Schrodinger}). If (\ref{small_F_0}) is satisfied, then for $a.e.$ $E\in \mathcal{I}$, the equation is reducible, i.e., there is a transformation $U(\omega t)$, unitary in $L^2$ and analytically depending on $t$, such that, by the transformation $u = U(\omega t) \, v$, Eq. (\ref{eq_Schrodinger}) is conjugated to the equation ${\rm i}\partial_t v = Gv$ with $G$ a linear operator independent of $t$. More precisely,
\begin{enumerate}
\item For a.e. $E\in\mathcal{I}\setminus\mathcal{O}_{\epsilon_0}$, G is unitary equivalent to $\frac{\varrho}{2}(x^2-\partial_x^2)$ for some $\varrho=\varrho_E\ge0$;
 \item If ${\rm Leb}(\Lambda_j) > 0$, then
 \begin{itemize}
 \item for $E\in {\rm int}\Lambda_j$,G is unitary equivalent to $-\frac{\lambda {\rm i}}{2}(x\cdot\partial_x+\partial_x\cdot x)$ for some $\lambda=\lambda_E>0$;
 \item for $E\in\partial\Lambda_j\setminus\partial\mathcal{I}$, G is unitary equivalent to $-\frac{\kappa}{2}x^2-{\rm i}\iota\partial_x$ for some $\kappa=\kappa_E\in\mathbb{R}\setminus\{0\}$ and $\iota=\iota_E\in\R$.
 \end{itemize}
If ${\rm Leb}(\Lambda_j)=0$, then $G$ is unitary equivalent to $-{\rm i}\iota\partial_x$
for some $\iota=\iota_E\in\R$ for $E\in\Lambda_j$.
\end{enumerate}
\end{Theorem}

\begin{remark} Compared with Theorem 6 of \cite{LZZ2021} on the reducibility of the quantum Hamiltonian with quadratic time quasi-periodic perturbation, the main difference in Theorem \ref{thm_Schro_reduc} lies in the linear term $-{\rm i}\iota\partial_x$.
\end{remark}

This remaining part of this section is devoted to the proof of Theorem \ref{thm_Schro_reduc}.

\subsection{Reducibility in classical hamiltonians}

Following the strategy of \cite{BGMR2018} and \cite{LZZ2021}, it is sufficient to consider the reducibility in the corresponding classical Hamiltonian systems for $E\in{\CI}$:
\begin{eqnarray}
h_E(\omega t, x,\xi)&=&\frac{\nu(E)}{2}(\xi^2+x^2)+W(E,\omega t, x,\xi)\label{classical_Hamiltonian}\\
&=&\frac12\left\la X, J (A_E+F_E(\omega t))X \right\ra + \left\la X, J b_E(\omega t)\right\ra,  \quad \nonumber
\end{eqnarray}
where, $X:=\left(\begin{array}{c}
                               x \\
                               \xi \\
                             \end{array}
                           \right)\in\R^2$, $J:=\left(\begin{array}{cc}
0 & -1 \\1 & 0 \end{array}\right)$, $A_E$ and $F_E$ are defined as in Section \ref{redu_qp_lin_sys}, and $b_E(\cdot):=\left(\begin{array}{c}
                               b_1(E,\cdot) \\[1mm]
                               b_2(E,\cdot)
                             \end{array}\right)\in C^\omega(\T^d,\R^2)$.
Since every quantity depends on $E\in{\CI}$, we do not always write this dependence explicitly in some statements in the remaining part of this paper. The following proposition is about the reducibility of the Hamiltonian $h_E$.

\begin{Proposition}\label{prop_reduc_ODE} {\bf (Reducibility between classical Hamiltonians)} If (\ref{small_F_0}) is satisfied, then, for a.e. $E\in{\CI}$, there are $K=K_E\in \N^*$ quadratic Hamiltonians
$$\chi_k(\omega t,x,\xi)= \frac12\left\la X, J Z_k(\omega t)X \right\ra, \quad Z_k\in C^\omega(2\T^d,{\rm sl}(2,\R)),\quad k=1,\cdots,K,$$
and one Hamiltonian
$$\chi_*(\omega t,x,\xi)=\la X, J l(\omega t) \ra + \varepsilon(\omega t),\qquad l\in C^\omega(2\T^d, \R^2), \quad \varepsilon \in C^\omega(2\T^d),$$
such that $h_E$ is conjugated to the Hamiltonian
$$g_E(x,\xi)= \frac12\left\la X, J B X \right\ra+ \la X, J w \ra+ {\CC},\qquad B\in {\rm sl}(2,\R),\quad w \in \R^2,\quad {\CC}\in\R,$$
via the composition of time$-1$ maps $\phi_{\chi_1}^1\circ\cdots\circ\phi_{\chi_K}^1\circ\phi_{\chi_*}^1$. More precisely,
\begin{enumerate}
\item For a.e. $E\in\mathcal{I}\setminus\mathcal{O}_{\epsilon_0}$, ${\rm det} B>0$ and $w=\left(\begin{array}{c}
0  \\[1mm] 0\end{array}\right)$.
 \item If ${\rm Leb}(\Lambda_j) > 0$, then
 \begin{itemize}
 \item for $E\in {\rm int}\Lambda_j$, ${\rm det} B<0$ and $w=\left(\begin{array}{c}
0  \\[1mm] 0\end{array}\right)$;
 \item for $E\in\partial\Lambda_j\setminus\partial\mathcal{I}$, $B$ is similar to
$\left(\begin{array}{cc}
0 & 0 \\[1mm] \kappa & 0 \end{array}\right)$ and $w=\left(\begin{array}{c}
\iota  \\[1mm] 0\end{array}\right)$
  for some $\kappa=\kappa_E\in\R\setminus\{0\}$ and $\iota=\iota_E\in\R$.
 \end{itemize}
If ${\rm Leb}(\Lambda_j)=0$, then for $E\in\Lambda_j$, $B=\left(\begin{array}{cc}
0 & 0 \\[1mm] 0 & 0 \end{array}\right)$ and $w=\left(\begin{array}{c}
\iota \\[1mm] 0\end{array}\right)$
  for some $\iota=\iota_E\in\R$.
\end{enumerate}
 \end{Proposition}

\begin{remark}
The linear term $-{\rm i}\iota\partial_x$ mentioned in Theorem \ref{thm_Schro_reduc} comes with $w=\left(\begin{array}{c}
\iota \\[1mm] 0\end{array}\right)$ after the conjugation in classical Hamiltonians if ${\rm det}B=0$.
The proof of Proposition \ref{prop_reduc_ODE} deals only the elimination of the time quasi-periodic linear terms for the classical Hamiltonian based on a constant quadratic part. As for the elimination of the time-dependent quadric terms, it can be found in Section 4 of \cite{LZZ2021}, (which is originally proposed by Eliasson \cite{Eli1992}) and we do not present it in details.

\end{remark}

\proof To show the reducibility of the classical Hamiltonian $h_E$, we can focus on the quasi-periodic affine system of equations of motion $(\omega, \, A+F(\cdot),  \, b(\cdot))$:
\begin{equation}\label{affine_sys}
X'(t)=(A+F(\omega t))X(t)+b(\omega t).
\end{equation}
In view of Proposition \ref{prop_eliasson}, we have the reducibility of the ${\rm sl}(2,\R)-$linear system $(\omega, \, A+F(\cdot))$ to the constant ${\rm sl}(2,\R)-$linear system $(\omega, \, B)$ by finitely many ${\rm SL}(2,\R)-$transformations, i.e., there exist $B\in {\rm sl}(2,\R)$ and $Z_j\in C^{\omega}(2\T^d,{\rm sl}(2,\R))$, $j=1,\cdots, K$, $K=K_E\in\N$, such that ${\CZ}:=\prod_{j=1}^K e^{Z_j}\in C^{\omega}(2\T^d,{\rm SL}(2,\R))$ satisfies that
\begin{equation}\label{reduci_Eliasson}
\frac{d}{dt}{\CZ}(\omega t)= (A+F(\omega t)){\CZ}(\omega t)-{\CZ}(\omega t) B,
\end{equation}
Under the transformation $X={\CZ}(\omega t) \tilde X$, the affine system (\ref{affine_sys}) is conjugated to
\begin{equation}\label{affine_to_reduce}
\tilde X'(t)=B\tilde X(t)+p(\omega t),\qquad   p(\omega t)=\left(\begin{array}{c}
p_1(\omega t) \\[1mm] p_2(\omega t)  \end{array}\right):={\CZ}(\omega t)^{-1}b(\omega t).
\end{equation}
As in \cite{LZZ2021},  we have the Hamiltonians
$$\chi_k(\omega t,x,\xi)= \frac12\left\la X, J Z_k(\omega t)X \right\ra,$$
which generates the time$-1$ maps $\phi_{\chi_k}^1(\omega t, x,\xi)$ such that $h_E$ is conjugated to
\begin{equation}
\tilde g(\omega t, x,\xi):=\frac12\left\la \tilde{X}, J B \tilde{X} \right\ra+\left\la \tilde{X}, J p(\omega t) \right\ra
\end{equation}
via the composition of maps $\phi_{\chi_1}^1(\omega t, x,\xi)\circ\cdots\circ\phi_{\chi_K}^1(\omega t, x,\xi)$. Here, we still use $(x,\xi)$ to denote the variables in $\R^2$ associated with $\tilde X$ for convenience.

We consider the reducibility of the affine system (\ref{affine_to_reduce}) via a transformation of the form
\begin{equation}\label{transform_affine}
\tilde X= R_\varphi Y+\left(
\begin{array}{l}
\alpha(\omega t)\\
\beta(\omega t)
\end{array}
\right),\end{equation}
with $R_\varphi\in {\rm SO}(2,\R)$ and $\alpha$, $\beta\in C^{\omega}(2\T^d)$ to be determined according to the type of constant matrix $B$ (hyperbolic, parabolic or elliptic). For convenience, we assume in the following that $B$ equals to its standard form.

\

\begin{itemize}
\item [{\bf Case 1.}] (hyperbolic case) $B=\left(\begin{array}{cc}
\lambda & 0 \\[1mm] 0 & -\lambda \end{array}\right)$ with $\lambda>0$
\end{itemize}

To find such $\alpha$ and $\beta$, it is sufficient to solve the equations
$$\frac{d}{dt}\alpha(\omega t)=p_1(\omega t)+\lambda\alpha(\omega t),\quad \frac{d}{dt}\beta(\omega t)=p_2(\omega t)-\lambda\beta(\omega t).$$
By the expansions in Fourier series of $p_1$, $p_2$, $\alpha$, $\beta$, i.e.,
$$p_1(\cdot)=\sum_{k\in\mathbb{Z}^d}\hat p_{1k}e^{\frac{\rm i}2\la k,\cdot\ra}, \ p_2(\cdot)=\sum_{k\in\mathbb{Z}^d}\hat p_{2k}e^{\frac{\rm i}2\la k,\cdot\ra}, \ \alpha(\cdot)=\sum_{k\in\mathbb{Z}^d}\hat \alpha_ke^{\frac{\rm i}2\la k,\cdot\ra}, \ \beta(\cdot)=\sum_{k\in\mathbb{Z}^d}\hat \beta_ke^{\frac{\rm i}2\la k,\cdot\ra},$$
the above two equations are equivalent to
$$\hat p_{1k}+\lambda\hat \alpha_k-\frac{\rm i}2\la k, \omega\ra \hat \alpha_k=0,\quad \hat p_{2k}-\lambda\hat \beta_k-\frac{\rm i}2\la k, \omega\ra \hat \beta_k=0,\qquad \forall \ k\in\Z^d.$$
Hence, by taking the Fourier coefficients
\begin{equation}\label{ab_hyp}
\hat \alpha_k=\frac{2\hat p_{1k}}{-2\lambda+{\rm i}\la k, \omega\ra}, \quad \hat \beta_k=\frac{2\hat p_{2k}}{2\lambda+{\rm i}\la k, \omega\ra},\qquad \forall \ k\in\Z^d,
\end{equation}
the affine system (\ref{affine_to_reduce}) is reduced to
$$Y'(t)=\left(\begin{array}{cc}
\lambda & 0 \\[1mm] 0 & -\lambda \end{array}\right) Y(t)$$
via the transformation of the form (\ref{transform_affine}) with $R_{\varphi}={\rm Id}$.
Since $\omega\in {\rm DC}_d(\gamma,\tau)$,
it is easy to see from (\ref{ab_hyp}) that $\alpha(\cdot),\beta(\cdot)\in C^\omega_{r'}(2\mathbb{T}^d)$ for any $0<r'<r$.

\

\begin{itemize}
\item [{\bf Case 2.}] (parabolic case) $B=\left(\begin{array}{cc}
0 & 0 \\[1mm] \kappa & 0 \end{array}\right)$ with $\kappa\in\R\setminus\{0\}$
\end{itemize}

It is sufficient to solve the equations
$$\frac{d}{dt}\alpha(\omega t)=p_1(\omega t)-\hat p_{10},\quad \frac{d}{dt}\beta(\omega t)=p_2(\omega t)+ \kappa\alpha(\omega t),$$
which are equivalent to the equations for Fourier coefficients
$$\kappa\hat \alpha_0+\hat p_{20}=0 \, ; \quad \frac{\rm i}2\la k,\omega\ra\hat \alpha_k-\hat p_{1k}=0,  \  \   \frac{\rm i}2\la k,\omega\ra\hat \beta_k-\hat p_{2k}-\kappa \hat \alpha_k=0,\quad \forall \ k\in\Z^d\setminus\{0\}.$$
Hence, by taking the Fourier coefficients
\begin{equation}\label{ab_para}
\hat \alpha_0=-\frac{\hat p_{20}}{\kappa} \, ; \quad \hat\alpha_k=\frac{2\hat p_{1k}}{{\rm i}\la k,\omega\ra},  \  \  \hat\beta_k=\frac{2\hat p_{2k}+2\kappa\hat\alpha_k}{{\rm i}\la k,\omega\ra}=-\frac{2{\rm i}\la k,\omega\ra \hat p_{2k}+4\kappa \hat p_{1k}}{\la k,\omega\ra^2},\quad \forall \ k\in\Z^d\setminus\{0\},
\end{equation}
the affine system (\ref{affine_to_reduce}) is reduced to
\begin{equation}\label{constant_stark}
Y'(t)=\left(\begin{array}{cc}
0 & 0 \\[1mm] \kappa & 0 \end{array}\right) Y(t)+\left(\begin{array}{c}
\iota \\[1mm] 0 \end{array}\right),\quad \iota:=\hat p_{10}
\end{equation}
via the transformation of the form (\ref{transform_affine}) with $R_{\varphi}={\rm Id}$.
Since $\omega\in {\rm DC}_d(\gamma,\tau)$,
it is easy to see from (\ref{ab_para}) that $\alpha(\cdot),\beta(\cdot)\in C^\omega_{r'}(2\mathbb{T}^d)$ for any $0<r'<r$.

\

\begin{itemize}
\item [{\bf Case 3.}] (elliptic case) $B= \left(\begin{array}{cc}
0 & \varrho \\[1mm] -\varrho & 0 \end{array}\right)$ with $\varrho\in\R$.
\end{itemize}

From Proposition \ref{prop_eliasson} and Remark \ref{remark_rot_num}, we have that $\varrho=\varrho_E$ satisfies the Diophantine condition w.r.t $\omega$, i.e. there exist ${\CK}={\CK}_E$ and $\sigma=\sigma_E>0$ such that
\begin{equation}\label{Dio_varrho}
\left|\varrho-\frac{\la n,\omega\ra}2\right|\ge\frac{\CK}{1+|n|^\sigma},\quad \forall \ n\in\mathbb{Z}^d \ {\rm  if} \ \varrho\neq 0,
\end{equation}
and if $\varrho=0$, (\ref{Dio_varrho}) holds obviously for every $n\in\mathbb{Z}^d\setminus\{0\}$.

If $\varrho\neq 0$, we try to solve the equations
$$\frac{d}{dt}\alpha(\omega t)=p_1(\omega t)+\varrho\beta(\omega t),\quad \frac{d}{dt}\beta(\omega t)=p_2(\omega t)-\varrho\alpha(\omega t),$$
which are equivalent to the equations for Fourier coefficients
$$\hat p_{1k}+\varrho\hat\beta_k-\frac{\rm i}2\la k,\omega\ra\hat\alpha_k=0, \quad \hat p_{2k}-\varrho\hat\alpha_k-\frac{\rm i}2\la k,\omega\ra\hat\beta_k=0, \qquad \forall \ k\in\Z^d.$$
Hence, by taking the Fourier coefficients
\begin{equation}\label{ab_elli-rho}
\hat\alpha_k=\frac{2{\rm i}\la k,\omega\ra \hat p_{1k}+4\varrho \hat p_{2k}}{4\varrho^2-\la k,\omega\ra^2},\quad \hat\beta_k=\frac{2{\rm i}\la k,\omega\ra \hat p_{2k}-4\varrho \hat p_{1k}}{4\varrho^2-\la k,\omega\ra^2}, \qquad  \forall \ k\in\Z^d,
\end{equation}
the affine system (\ref{affine_to_reduce}) is reduced to
$$Y'(t)=\left(\begin{array}{cc}
0 & \varrho \\[1mm] -\varrho & 0 \end{array}\right) Y(t)$$
via the transformation of the form (\ref{transform_affine}) with $R_{\varphi}={\rm Id}$.
By the Diophantine condition (\ref{Dio_varrho}) of $\varrho$, we have, for $ k\in\Z^d$,
\begin{eqnarray}
|\hat\alpha_k|&\le&\frac{2|\la k,\omega\ra|}{|4\varrho^2-\la k,\omega\ra^2|}\left|\hat p_{1k}\right|+\frac{4|\varrho|}{|4\varrho^2-\la k,\omega\ra^2|}\left|\hat p_{2k}\right|\label{esti_alpha_k}\\
&\le& \frac{1+ |k|^\sigma}{2\CK}\left(\left|\hat p_{1k}\right|+\left|\hat p_{2k}\right|\right)\nonumber\\
&\le& \frac{1+ |k|^\sigma}{2\CK}e^{-\frac{|k|r}{2}}\left(|p_1|_r+|p_2|_r\right),\nonumber\\
|\hat\beta_k|&\le&\frac{1+ |k|^\sigma}{2\CK}e^{-\frac{|k|r}{2}}\left(|p_1|_r+|p_2|_r\right).\label{esti_beta_k}\end{eqnarray}

If $\varrho=0$, we try to solve the equations
$$\frac{d}{dt}\alpha(\omega t)=p_1(\omega t)-\hat p_{10}, \quad \frac{d}{dt}\beta(\omega t)=p_2(\omega t)-\hat p_{20},$$
which are equivalent to the equations for Fourier coefficients
$$\hat p_{1k}-\frac{\rm i}2\la k,\omega\ra\hat{\alpha}_k=0, \quad \hat p_{2k}-\frac{\rm i}2\la k,\omega\ra\hat{\beta}_k=0, \qquad \forall \ k\in\Z^d\setminus\{0\}.$$
By taking the Fourier coefficients
\begin{equation}\label{alpha_beta_k-0}
\hat{\alpha}_k=-\frac{2{\rm i} \,  \hat p_{1k} }{\la k,\omega\ra},\quad \hat{\beta}_k=-\frac{2{\rm i}  \, \hat p_{2k} }{\la k,\omega\ra},
\qquad  \forall \ k\in\Z^d\setminus\{0\},
\end{equation}
the affine system (\ref{affine_to_reduce}) is reduced to
$$ \tilde Y'(t)=\left(\begin{array}{c}
\hat p_{10} \\[1mm] \hat p_{20} \end{array}\right),$$
via the transformation $$\tilde X= \tilde Y+\left(
\begin{array}{l}
\alpha(\omega t)\\
\beta(\omega t)
\end{array}
\right).$$
Then, by the rotation
$$\tilde Y =\frac{1}{\sqrt{\hat p_{10}^2+\hat p_{20}^2}}\left(\begin{array}{cc}
\hat p_{10} & -\hat p_{20} \\[1mm] \hat p_{20} & \hat p_{10} \end{array}\right)Y,$$
we obtain the system
\begin{equation}\label{constant_degenerate}
Y'(t)=\left(\begin{array}{c}
\iota \\[1mm] 0 \end{array}\right),\quad \iota:=\sqrt{\hat p_{10}^2+\hat p_{20}^2}.
\end{equation}
Hence, (\ref{affine_to_reduce}) is reduced to (\ref{constant_degenerate}) via the transformation of the form (\ref{transform_affine}) with
\begin{equation}\label{R_phi}
R_{\varphi}:=\frac{1}{\sqrt{\hat p_{10}^2+\hat p_{20}^2}}\left(\begin{array}{cc}
\hat p_{10} & -\hat p_{20} \\[1mm] \hat p_{20} & \hat p_{10} \end{array}\right).
\end{equation}

By (\ref{esti_alpha_k}) -- (\ref{alpha_beta_k-0}), we see that, for all $\varrho\in\R$, $\alpha(\cdot)$, $\beta(\cdot)\in C^\omega_{r'}(2\mathbb{T}^d)$ for any $0<r'<r$.

\

Combining all above cases, we get the reducibility of the affine system (\ref{affine_sys}) via the transformation of the form (\ref{transform_affine}). It is summarized in the following table.
{\center{\begin{tabular}{|l|c|c|}
\hline
&  $R_{\varphi}$ & coefficients of $\alpha(\cdot)$, $\beta(\cdot)$  \\[2mm]
 \hline
$B=\left(\begin{array}{cc}
\lambda & 0 \\[1mm] 0 & -\lambda \end{array}\right)$ with $\lambda>0$  & {\rm Id} & given in (\ref{ab_hyp})  \\[2mm]
\hline
$B=\left(\begin{array}{cc}
0 & 0 \\[1mm] \kappa & 0 \end{array}\right)$ with $\kappa\in\R\setminus\{0\}$ & {\rm Id} & given in (\ref{ab_para})  \\[2mm]
\hline
$B=\left(\begin{array}{cc}
0 & \varrho \\[1mm] -\varrho & 0 \end{array}\right) $ with $\varrho$ satisfying (\ref{Dio_varrho}) & {\rm Id} & given in (\ref{ab_elli-rho})  \\[2mm]
\hline
$B=\left(\begin{array}{cc}
0 & 0 \\[1mm] 0 & 0 \end{array}\right) $  & given in (\ref{R_phi}) & given in (\ref{alpha_beta_k-0}) \\[2mm]
\hline
\end{tabular}}}\\
Since the rotation $R_\varphi$ satisfies that
$$R_\varphi=\exp\left\{\left(\begin{array}{cc}
0 & -\varphi \\[1mm] \varphi & 0\end{array}\right)\right\},$$
it can be seen as a constant member of $\left\{e^{Z_j(\cdot)}\right\}_{j=1}^K$.

\begin{remark} We can see from (\ref{constant_stark}) and (\ref{constant_degenerate}) that the constant $\iota$ comes from the average term of  $p(\cdot)={\CZ}(\cdot)^{-1}b(\cdot)$. The transformation ${\mathcal Z}(\cdot)$ is exactly that constructed by the scheme of Eliasson \cite{Eli1992} and usually we do not have its precise form, hence it is not possible to compute $\iota$. However, as mentioned in Remark \ref{rem_iota}, $\iota$ is generically non-vanishing. Indeed, a perturbation on $b(\cdot)$ will easily cause a non-vanishing average term and hence a non-vanishing $\iota$.
\end{remark}

Now it remains to translate the conjugation between affine systems into that of classical Hamiltonians. More generally, we consider two classical Hamiltonians
$$h_j(\omega t, x,\xi)=\frac12\left\la X, J A(\omega t)X \right\ra +\left\la X, J b_j(\omega t) \right\ra,  \quad j=1,2$$
with $X=\left(\begin{array}{c}
                               x \\[1mm]
                               \xi
                             \end{array}\right)$, $A(\cdot)\in C_r^{\omega}(\T^d, {\rm sl}(2,\R))$ and $b_j(\cdot)\in C_r^\omega(2\T^d,\R^2)$.
The corresponding equations of motion are
$$x'=\frac{\partial h_j}{\partial\xi},\quad \xi'=-\frac{\partial h_j}{\partial x},\qquad j=1,2,$$
which are the affine systems $(\omega, \, A(\cdot), \, b_j(\cdot))$:
$$\left(\begin{array}{c}
          x(t) \\
          \xi(t)
        \end{array}
\right)'=A(\omega t)\left(\begin{array}{c}
          x(t) \\
          \xi(t)
        \end{array}
\right)+b_j(\omega t).$$
The following lemma completes the proof of Proposition \ref{prop_reduc_ODE}.

\begin{Lemma}\label{lemma_ham_cl-1}
If the affine system $\big(\omega, \, A(\cdot), \, b_1(\cdot)\big)$ is conjugated to
$\big(\omega, \, A(\cdot), \, b_2(\cdot)\big)$ by a time quasi-periodic $\R^2-$transformation $l(\omega t)$, i.e.,
\begin{equation}\label{conj_ode-1}
\frac{d}{dt}l(\omega t)=A(\omega t)l(\omega t)+b_1(\omega t)-b_2(\omega t),\qquad l \in C^\omega(2\T^d, \R^2),
\end{equation}
then the classical Hamiltonian $h_1(\omega t,x,\xi)$ is conjugated to $h_2(\omega t,x,\xi)+{\CC}$ via the time$-1$ flow $\phi_{\chi_*}^1(t,x,\xi)$
 generated by the Hamiltonian
 \begin{equation}\label{chi_l}
 \chi_*(\omega t,x,\xi)=\la X, Jl(\omega t)\ra +\varepsilon(\omega t),
 \end{equation}
 where, with $F(\theta) : =\frac12  \la l(\theta), J (b_1(\theta)+b_2(\theta))\ra=  \sum\limits_{k\in \Z^d} \hat{F}_k e^{\frac{\rm i}{2}\la k, \theta\ra }$, we have
 $${\CC}:= \hat{F}_0 = \frac{1}{(4\pi)^d}\int_{2\mathbb{T}^d} F(\theta)  \, d\theta,\quad \varepsilon(\theta):=-2{\rm i}\sum_{k\in\Z^d\setminus\{0\}}\frac{\hat{F}_k}{ \la k, \omega\ra } e^{\frac{\rm i}{2}\la k, \theta\ra}.$$
\end{Lemma}
\begin{remark}
It is obvious that the constant ${\CC}$ in the classical Hamiltonian does not influence the growth of Sobolev norm in the quantized Hamiltonian. Hence we usually say that $h_1$ is conjugated to $h_2$ by ignoring this constant.
\end{remark}
\proof The time$-\tau$ flow generated by $\chi_*(\omega t,x,\xi)$ in (\ref{chi_l}) is
$\phi_{\chi_*}^\tau(t,x,\xi)=X+\tau l(\omega t)$.
In view of (\ref{time1}), we see that $h_1(\omega t,x,\xi)$ is conjugated to
\begin{eqnarray*}
 & &(h_1\circ\phi^1_{\chi_{*}})(\omega t, x,\xi)-\int_0^1\frac{\partial{\chi_*}}{\partial t}\left(\omega t,\phi_{\chi_*}^\tau(t,x,\xi)\right) \, d\tau\\
& =& \frac12\la X+l(\omega t),JA(\omega t)(X+l(\omega t))\ra +\la X+l(\omega t), J b_1(\omega t)\ra \\
& & - \, \int_0^1 \left(\left\la X+\tau l(\omega t),J \frac{d}{dt}l(\omega t)\right\ra+\frac{d}{dt} \varepsilon(\omega t) \right) \, d\tau\\
&= &\frac12\la X,JA(\omega t)X\ra + \la X,J b_2(\omega t)\ra\\
& & + \, \frac12\la l(\omega t),J(b_1(\omega t)+b_2(\omega t))\ra-\frac{d}{dt} \varepsilon(\omega t)\\
&= &h_2(\omega t, x, \xi) + {\CC},
 \end{eqnarray*}
 by noting that $JA(\omega t)$ is symmetric when $A(\omega t)\in {\rm sl}(2,\R)$.\qed

\subsection{Reducibility in quantum hamiltonians -- Proof of Theorem \ref{thm_Schro_reduc}}\label{sec_proof_reduc}

Note that the time-dependent Schr\"odinger equation (\ref{eq_Schrodinger}) is
$${\rm i}\partial_t u=h_E^W(\omega t ,x, -{\rm i}\partial_x) u.$$
According to Proposition \ref{Prop_hami} and \ref{prop_reduc_ODE},
if (\ref{small_F_0}) is satisfied, then, by the unitary transformation
$$u=U(\omega t)  \,  v,\quad  U(\omega t):=e^{-{\rm i}\chi_1^W(\omega t,x,\xi)}\cdots e^{-{\rm i}\chi_K^W(\omega t,x,\xi)}e^{-{\rm i}\chi_*^W(\omega t,x,\xi)},$$
Eq. (\ref{eq_Schrodinger}) is conjugated to the equation
\begin{equation}\label{reduced_quantum_ham}
{\rm i}\partial_t v = G(x, -{\rm i}\partial_x) \, v, \qquad G(x, -{\rm i}\partial_x):=g^W(x, -{\rm i}\partial_x),
\end{equation}
with the symbol of $G$ of the form
$$g(x,\xi)= \frac12\left\la X, J B X \right\ra+ \la X, J w \ra+ {\CC},\qquad B\in {\rm sl}(2,\R),\quad w \in \R^2,\quad {\CC}\in\R.$$

Let us focus on the concrete form of $G$ for different parameter $E$. Since the similarity between ${\rm sl}(2,\R)$ matrices is related to the unitary equivalence between quantized Hamiltonians, we have the following conclusions.
\begin{enumerate}
\item For a.e. $E\in\mathcal{I}\setminus\mathcal{O}_{\epsilon_0}$, $B$ is similar to
$\left(\begin{array}{cc}
            0 & \varrho\\[1mm]
            -\varrho & 0
          \end{array}
\right)$ for some $\varrho=\varrho_E\geq 0$ and $w=\left(\begin{array}{c}
0  \\[1mm] 0\end{array}\right)$, hence $G$ is unitary equivalent to the Weyl quantization
of $$\frac12\left\la X, J \left(\begin{array}{cc}
            0 & \varrho\\[1mm]
            -\varrho & 0
          \end{array}
\right) X \right\ra=\frac{\varrho}2(x^2+\xi^2),$$
which is $\frac{\varrho}{2}(x^2-\partial_x^2)$.\\[1mm]

 \item If ${\rm Leb}(\Lambda_j) > 0$, then
 \begin{itemize}
 \item for $E\in {\rm int}\Lambda_j$, $B$ is similar to $\left(\begin{array}{cc}
            \lambda & 0\\[1mm]
            0 & -\lambda
          \end{array}
\right)$ for some $\lambda=\lambda_E>0$ and $w=\left(\begin{array}{c}
0  \\[1mm] 0\end{array}\right)$, hence $G$ is unitary equivalent to
$-\frac{\lambda {\rm i}}{2}(x\cdot\partial_x+\partial_x\cdot x)$.\\[1mm]

 \item for $E\in\partial\Lambda_j\setminus\partial\mathcal{I}$, $B$ is similar to
$\left(\begin{array}{cc}
0 & 0 \\[1mm] \kappa & 0 \end{array}\right)$ and $w=\left(\begin{array}{c}
\iota  \\[1mm] 0\end{array}\right)$
  for some $\kappa=\kappa_E\in\R\setminus\{0\}$ and $\iota=\iota_E\in\R$, hence $G$ is unitary equivalent to the Stark Hamiltonian
  \begin{equation}\label{stark}
  -\frac{\kappa}{2}x^2-{\rm i}\iota\partial_x.
  \end{equation}
 \end{itemize}
If ${\rm Leb}(\Lambda_j)=0$, then for $E\in\Lambda_j$, $B=\left(\begin{array}{cc}
0 & 0 \\[1mm] 0 & 0 \end{array}\right)$ and $w=\left(\begin{array}{c}
\iota \\[1mm] 0\end{array}\right)$
  for some $\iota=\iota_E\in\R$, hence $G$ is unitary equivalent to $-{\rm i}\iota\partial_x$.
\end{enumerate}

\section{Growth of Sobolev norms via reducibility}\label{sec_Growth}


According to (\ref{norm_U}), to precise the growth of Sobolev norms for the solution to Eq. (\ref{eq_Schrodinger}), it is sufficient to study the reduced quantum Hamiltonian $G(x, -{\rm i}\partial_x)$ obtained in Theorem \ref{thm_Schro_reduc}, or more simply, the unitary equivalent forms listed in Theorem \ref{thm_Schro_reduc}.

\smallskip

For a.e. $E\in\mathcal{I}\setminus\mathcal{O}_{\epsilon_0}$, $G$ is unitary equivalent to $\frac{\varrho}{2} H_0$. Since the ${\CH}^s-$norm of $e^{-{\rm i}t\frac{\varrho}{2} H_0}\psi_0$ is conserved for any $\psi_0\in{\CH}^s$, the boundedness of Sobolev norm is shown.
Now we focus on the other cases, in which the growth of Sobolev norm occurs.

\begin{Proposition} \label{prop6}
For the equation
$$\partial_t v(t,x)=-\frac\lambda2 x\cdot\partial_x v (t,x)-\frac\lambda2 \partial_x(x\cdot v (t,x)), \qquad \lambda>0,$$
 with non-vanishing initial condition $ v (0, \cdot)= v_0(\cdot)\in {\CH}^s$, $s\geq 0$, there exist two constants $\tilde c, \, \tilde C>0$, depending on $s$, $\lambda$ and $ v _0$, such that the solution satisfies
$$\tilde c e^{\lambda st} \leq \|v(t,\cdot )\|_{s}\leq   \tilde C e^{\lambda st}, \qquad \forall \ t\geq0.$$
\end{Proposition}
The above proposition gives the exponential growth of Sobolev norm for $E\in{\rm int}\Lambda_j$ when ${\rm Leb}(\Lambda_j)>0$.
For the proof, see Proposition 2 of \cite{LZZ2021} or Remark 1.4 of \cite{MR2017}.



\

\begin{Proposition}\label{prop_para}
For the equation
\begin{equation}\label{eq_para}
{\rm i}\partial_t v (t,x)=-\frac{\kappa}{2} x^2\cdot v (t,x)-{\rm i}\iota\partial_x v(t,x), \qquad \kappa,\iota\in\R\setminus\{0\},
\end{equation}
with non-vanishing initial condition $v (0, \cdot)= v_0(\cdot)\in {\CH}^s$, $s\geq 0$, there exists constants $\tilde c, \, \tilde C>0$, depending on $s$, $\kappa$, $\iota$ and $v _0$, such that the solution satisfies
\begin{equation}\label{sobo_para}
\tilde c |\iota\kappa|^s t^{2s} \leq \|v(t,\cdot)\|_{s}\leq \tilde C \left(1+|\iota\kappa|^s t^{2s}\right), \quad \forall \ t\geq 0.
\end{equation}
\end{Proposition}

The above proposition gives the possible $t^{2s}-$growth of ${\CH}^s$ norm of solutions for $E \in \partial \Lambda_j\setminus\partial{\CI}$ when ${\rm Leb}(\Lambda_j) > 0$. It is indispensable for the $t^{2s}-$growth that the coefficient $\iota$ in the reduced constant quantum Hamiltonian is non-vanishing.
For the case $\iota=0$, the $\mathcal{H}^s-$norm of the solution to Eq. (\ref{eq_para}) presents $t^s-$growth (see Proposition 3 of \cite{LZZ2021}).

\proof
With the initial condition $v(0,\cdot)=v_0(\cdot)\in\mathcal{H}^s$, it can be verified that the solution to Eq. (\ref{eq_para}) is
$$v(t,x)=e^{\frac{{\rm i}\kappa}{6\iota}x^3}e^{-\frac{{\rm i}\kappa}{6\iota}(x-\iota t)^3}v_0(x-\iota t).$$
For any $s\ge0$, we have
$$\|x^s v(t,x)\|_{L^2}^2=\int_\mathbb{R}x^{2s}|v(t,x)|^2dx
=\int_\mathbb{R}x^{2s}|v_0(x-\iota t)|^2dx
=\int_\mathbb{R}(y+\iota t)^{2s}|v_0(y)|^2dy.$$
Hence, there exist two constants $c'$, $C'>0$, depending on $s$, $\iota$ and $v_0$, such that
\begin{equation}\label{growth_part1}
c' (1+|\iota|^s t^{s})\leq \|x^sv(t,x)\|_{L^2} \leq C'(1+|\iota|^s t^{s}) .
\end{equation}


In view of the equivalent definitions of the ${\CH}^s-$norm given in (\ref{norm_equiv-1}), it remains to describe the growth of $\|\partial_x^s v(t,x)\|_{L^2}$. We first consider the derivatives of
$e^{\frac{{\rm i}\kappa}{6\iota}\left(3y^2\iota t+3y\iota^2t^2+\iota^3t^3\right)}$ w.r.t. $y\in\R$. We have that, for $\alpha\in\mathbb{N}^*$, \begin{equation}\label{dev_poly}
\partial_y^\alpha e^{\frac{{\rm i}\kappa}{6\iota}\left(3y^2\iota t+3y\iota^2t^2+\iota^3t^3\right)}=P_{2\alpha}(t,y) \, e^{\frac{{\rm i}\kappa}{6\iota}\left(3y^2\iota t+3y\iota^2t^2+\iota^3t^3\right)},
\end{equation}
where $P_{2\alpha}(t,y)$ is a multivariate polynomial of $(t,y)$ of degree $2\alpha$, satisfying
\begin{equation}\label{P_2alpha}
P_{2\alpha}(t,y)=\left(\frac{\rm i}{2}\iota\kappa \right)^\alpha t^{2\alpha}+\sum_{1 \le a\le 2\alpha-1}Q_{\alpha,a}(y) \, t^a,
\end{equation}
with $Q_{\alpha,a}(y)$ a polynomial of $y$ of degree no more than $2\alpha-1$.
Indeed, for $\alpha=1,2$, direct computations show that
\begin{eqnarray*}
\partial_y e^{\frac{{\rm i}\kappa}{6\iota}(3y^2\iota t+3y\iota^2t^2+\iota^3t^3)}&=&\left(\frac{\rm i}{2}\iota\kappa t^2+{\rm i}\kappa yt\right) e^{\frac{{\rm i}\kappa}{6\iota}(3y^2\iota t+3y\iota^2t^2+\iota^3t^3)},\\
\partial^2_y e^{\frac{{\rm i}\kappa}{6\iota}(3y^2\iota t+3y\iota^2t^2+\iota^3t^3)}&=&\left(\left(\frac{\rm i}{2}\iota\kappa t^2\right)^2-\iota \kappa^2y t^3-\kappa^2y^2 t^2+{\rm i}\kappa t\right) e^{\frac{{\rm i}\kappa}{6\iota}(3y^2\iota t+3y\iota^2t^2+\iota^3t^3)},\\
\end{eqnarray*}
and if (\ref{dev_poly}) holds for some $\alpha\in\mathbb{N}^*$, $\alpha\geq 2$, with $P_{2\alpha}$ satisfying (\ref{P_2alpha}), then
\begin{eqnarray*}
& & \partial_y^{\alpha+1} e^{\frac{{\rm i}\kappa}{6\iota}(3y^2\iota t+3y\iota^2t^2+\iota^3t^3)}\\
&=&\partial_y\left((P_{2\alpha}(t,y) \, e^{\frac{{\rm i}\kappa}{6\iota}(3y^2\iota t+3y\iota^2t^2+\iota^3t^3)}\right)\\
&=&\left( \left(\frac{\rm i}2\iota \kappa t^2+{\rm i}\kappa yt\right)P_{2\alpha}(t,y)+\partial_yP_{2\alpha}(t,y)\right)e^{\frac{{\rm i}\kappa}{6\iota}(3y^2\iota t+3y\iota^2t^2+\iota^3t^3)}\\
&=&\left(\left( \frac{\rm i}2\iota \kappa \right)^{\alpha+1}t^{2\alpha+2}+{\rm i}\kappa \left( \frac{\rm i}2\iota \kappa \right)^{\alpha} y t^{2\alpha+1}+ \left(\frac{\rm i}2\iota \kappa t^2+{\rm i}\kappa yt\right) \sum_{1\le a\le 2\alpha-1}Q_{\alpha,a}(y) \, t^a  \right.
\\
& & \  \  \  \  \  \   \left. + \, \sum_{1\le a\le 2\alpha-1}Q'_{\alpha,a}(y) \, t^a \right)   \, e^{\frac{{\rm i}\kappa}{6\iota}(3y^2\iota t+3y\iota^2t^2+\iota^3t^3)} \\
&=&\left(\left( \frac{\rm i}2\iota \kappa \right)^{\alpha+1}t^{2\alpha+2}+  \sum_{1\le a\le 2\alpha+1} Q_{\alpha+1,a}(y) \, t^a
\right) \, e^{\frac{{\rm i}\kappa}{6\iota}(3y^2\iota t+3y\iota^2t^2+\iota^3t^3)}\\
&=&P_{2\alpha+2}(t,y) \, e^{\frac{{\rm i}\kappa}{6\iota}(3y^2\iota t+3y\iota^2t^2+\iota^3t^3)},
\end{eqnarray*}
with the polynomials $Q_{\alpha+1,a}(y)$ defined as
$$Q_{\alpha+1,a}(y)=\left\{\begin{array}{cc}
{\rm i}\kappa \left( \frac{\rm i}2\iota \kappa \right)^{\alpha} y+ \frac{\rm i}2\iota \kappa Q_{\alpha,2\alpha-1}(y), & a= 2\alpha+1 \\[1mm]
\frac{\rm i}2\iota \kappa Q_{\alpha,2\alpha-2} (y) +  {\rm i}\kappa y Q_{\alpha,2\alpha-1}(y), & a= 2\alpha \\[1mm]
\frac{\rm i}2\iota \kappa Q_{\alpha,a-2}(y)+  {\rm i}\kappa y Q_{\alpha,a-1}(y) + Q'_{\alpha,a}(y), & 3 \leq a\leq 2\alpha-1
\\[1mm]
  {\rm i}\kappa y Q_{\alpha,1}(y) + Q'_{\alpha,2}(y), &  a=2
\\[1mm]
Q'_{\alpha,1}(y),& a=1
\end{array} \right. .$$
Note that, for $y=x-\iota t$,
$$\partial_x^s v(t,x)=\partial_x^s \left(e^{\frac{{\rm i}\kappa}{6\iota}x^3}e^{-\frac{{\rm i}\kappa}{6\iota}(x-\iota t)^3}v_0(x-\iota t)\right) =\partial_y^s \left(e^{\frac{{\rm i}\kappa}{6\iota}(3y^2\iota t+3y\iota^2t^2+\iota^3t^3)}v_0(y)\right).$$
Then, by (\ref{dev_poly}) and (\ref{P_2alpha}), we have, for $s\in\mathbb{N}^*$,
\begin{eqnarray*}
& &\left|\partial_x^s v(t,x)-\left(\frac{{\rm i}\iota\kappa}{2}\right)^st^{2s} e^{\frac{{\rm i}\kappa}{6\iota}(3y^2\iota t+3y\iota^2t^2+\iota^3t^3)}v_0(y)\right|\\
&=&\left|\sum_{0\le\alpha\le s}C_s^\alpha\left(\partial_y^\alpha e^{\frac{{\rm i}\kappa}{6\iota}(3y^2\iota t+3y\iota^2t^2+\iota^3t^3)}\right) \left(\partial_y^{s-\alpha}v_0(y)\right)-\left(\frac{{\rm i}\iota\kappa}{2}\right)^st^{2s}e^{\frac{{\rm i}\kappa}{6\iota}(3y^2\iota t+3y\iota^2t^2+\iota^3t^3)}v_0(y)\right|\\
&=&\left|\sum_{0\le\alpha\le s}C_s^\alpha P_{2\alpha}(t,y) \, \partial_y^{s-\alpha}v_0(y)-\left(\frac{{\rm i}\iota\kappa}{2}\right)^st^{2s} v_0(y)\right|\\
&=&\left|\sum_{0\le\alpha\le s-1 }C_s^\alpha P_{2\alpha}(t,y) \, \partial_y^{s-\alpha}v_0(y)+v_0(y)\sum_{1 \le a\le 2s-1}Q_{s,a}(y) \, t^a \right|
\end{eqnarray*}
Hence, there exists a constant $c''$, depending on $s$, $\iota$, $\kappa$ and $v_0$, such that
$$ \left\|\partial_x^s v(t,x)-\left(\frac{{\rm i}\iota\kappa}{2}\right)^st^{2s} e^{\frac{{\rm i}\kappa}{6\iota}(3y^2\iota t+3y\iota^2t^2+\iota^3t^3)}v_0(y)\right\|_{L^2} \leq c'' \left(1+t^{2s-1}\right).$$
On the other hand, we have
$$\left\|\left(\frac{{\rm i}\iota\kappa}{2}\right)^st^{2s} e^{\frac{{\rm i}\kappa}{6\iota}(3y^2\iota t+3y\iota^2t^2+\iota^3t^3)}v_0(y)\right\|_{L^2}=\left(\frac{|\iota\kappa|}{2}\right)^st^{2s}\|v_0\|_{L^2},$$
Then, as $t\to\infty$,
\begin{equation}\label{growth_part2}
\|\partial_x^s v(t,x)\|_{L^2}\sim \frac{|\iota\kappa|^s}{2^s}t^{2s}\|v_0\|_{L^2}.
\end{equation}

Combining (\ref{growth_part1}) and (\ref{growth_part2}), we get the growth of ${\CH}^s-$norm in (\ref{sobo_para}). \qed

\

As for the degenerate case $E\in \Lambda_j$ with ${\rm Leb}(\Lambda_j)=0$, the reduced quantum Hamiltonian gives the equation
$$\partial_t v(t,x)=-\iota \partial_x v(t,x).$$
It is proved that the solution is the traveling wave
$$v(t,x)=v_0(x-\iota t),\quad v_0 \in {\CH}^s,$$
and its ${\CH}^s-$norm presents $t^s-$growth (see Proof of Corollary A.2 of \cite{BGMR2018}). This completes the proof of Theorem \ref{thm_Schro_sobolev}.

\

\noindent
{\bf Proof of Theorem \ref{super-ballistic}.} By Fourier transform
$$u(x)\mapsto  \hat u(\xi):=  \int_\R  e^{-2\pi {\rm i} \xi x} u(x) \, dx,$$
the linear equation
$${\rm i}\partial_t u(t,x)=({\CS}_au)(t,x)=-\partial_x^2 u(t,x)+ a x u(t,x),\quad u(0,x)=u_0(x),\quad x\in\R,$$
can be transformed into
$${\rm i}\partial_t \hat u(t,\xi)=4\pi^2 \xi^2 \hat u(t,\xi)+ {\rm i} \frac{a}{2\pi}  \partial_{\xi} \hat u(t,\xi),\quad \hat u(0,\xi)=\hat u_0(\xi),\quad \xi \in\R.$$
In view of (\ref{growth_part2}) in the proof of Proposition \ref{prop_para}, we have
$$\|x^s u(t,x)\|_{L^2(dx)}=\frac{1}{(2\pi)^s}\|\partial_\xi^s \hat u(t,\xi)\|_{L^2(d\xi)}\sim |a|^s t^{2s}\|u_0\|_{L^2(dx)},  \quad t\rightarrow \infty. \qed$$

\section{Proof of Theorem \ref{examthm}}\label{sec_pr_examples}

We show the growth of Sobolev norms of solutions to Eq. (\ref{exameq1}) by its reducibility.

\smallskip

Through Weyl quantization, the reducibility of Eq. (\ref{exameq1}) is equivalent to the reducibility of the corresponding classical Hamiltonian system:
\begin{align*}
\frac{d}{dt}\left(
\begin{array}{l}
x\\
\xi
\end{array}
\right)&=\left(\left(
\begin{matrix}
0 & 1\\
-1 & 0
\end{matrix}
\right)+\kappa\left(
\begin{matrix}
\cos(t)\sin(t) & -\sin^2(t) \\
\cos^2(t)  & -\cos(t) \sin(t)
\end{matrix}
\right)\right)\left(
\begin{array}{l}
x\\
\xi
\end{array}
\right)-2\iota\left(
\begin{array}{cc}
\cos(t)\\
0
\end{array}
\right)\\
&=\left(-J+\kappa e^{-t J}\left(
\begin{matrix}
0 &0\\
\cos(t)  & -\sin(t)
\end{matrix}
\right)\right)\left(
\begin{array}{l}
x\\
\xi
\end{array}
\right)-2\iota\left(
\begin{array}{cc}
\cos(t) \\
0
\end{array}
\right),
\end{align*}
where $J=\left(\begin{matrix}0 & -1\\1 & 0\end{matrix}\right).$
Under the transformation
$$\left(
\begin{array}{l}
x\\
\xi
\end{array}
\right)=e^{-t J}\left(
\begin{array}{l}
\tilde{x}\\
\tilde{\xi}
\end{array}
\right),$$
the above system is conjugated to
\begin{align*}
\frac{d}{dt}\left(
\begin{array}{l}
\tilde{x}\\
\tilde{\xi}
\end{array}
\right)&=e^{t J}\left(-J+\kappa e^{-t J}\left(
\begin{matrix}
0 &0\\
\cos(t) & -\sin(t)
\end{matrix}
\right)+J\right)e^{-t J}\left(
\begin{array}{l}
\tilde{x}\\
\tilde{\xi}
\end{array}
\right)-2\iota e^{t J}\left(
\begin{array}{cc}
\cos(t)\\
0
\end{array}
\right)\\
&=\kappa\left(
\begin{matrix}
0 &0\\
\cos(t) & -\sin(t)
\end{matrix}
\right)e^{-t J}\left(
\begin{array}{l}
\tilde{x}\\[1mm]
\tilde{\xi}
\end{array}
\right)-2\iota e^{t J}\left(
\begin{array}{cc}
\cos(t)\\[1mm]
0
\end{array}
\right)\\
&=\left(
\begin{matrix}
0 & 0\\[1mm]
\kappa & 0
\end{matrix}
\right)\left(
\begin{array}{l}
\tilde{x}\\[1mm]
\tilde{\xi}
\end{array}
\right)-2\iota\left(
\begin{array}{cc}
\cos^2(t)\\[1mm]
\sin(t)\cos(t)
\end{array}
\right).
\end{align*}
Then, by the transformation
$$\left(
\begin{array}{l}
\tilde{x}\\[1mm]
\tilde{\xi}
\end{array}
\right)=\left(
\begin{array}{l}
x_\ast\\[1mm]
\xi_\ast
\end{array}
\right)-\left(
\begin{array}{l}
\iota\sin(t)\cos(t)\\[1mm]
\left(\frac12\kappa\iota+\iota\right)\sin^2(t)
\end{array}
\right),$$
we get the affine system
\begin{equation}\label{eq2}
\frac{d}{dt}\left(
\begin{array}{l}
x_\ast\\
\xi_\ast
\end{array}
\right)=\left(
\begin{matrix}
0 & 0\\
\kappa & 0
\end{matrix}
\right)\left(
\begin{array}{l}
x_\ast\\
\xi_\ast
\end{array}
\right)-\left(
\begin{array}{l}
\iota\\
0
\end{array}
\right).
\end{equation}

An explicit computation shows that the two above transformations are respectively the time $1-$flows of Hamiltonians
\begin{eqnarray*}
\chi_1(t,x,\xi)&=&\frac{t}2(x^2+\xi^2),\\
\chi_2(t,x,\xi)&=&\left(\frac12\kappa\iota+\iota\right)\sin^2(t)\cdot x-\iota\sin(t)\cos(t)\cdot\xi.
\end{eqnarray*}
Hence, by the unitary transformation
$$u =e^{-{\rm i}\chi_1^W(t,x,-{\rm i}\partial_x)}e^{-{\rm i}\chi_2^W(t,x,-{\rm i}\partial_x)} v,$$
Eq. (\ref{exameq1}) is reduced to
$$
{\rm i}\partial_tv(t,x)=-\frac\kappa2x^2 v(t,x)-{\rm i}\iota\partial_x v(t,x).
$$
Then Theorem \ref{examthm} follows from Proposition \ref{prop_para}.

\appendix


\end{document}